\documentclass[11pt,leqno]{article}
\title{\ }
\author{\ }
\begin{document}
\newtheorem{theorem}{Theorem}[section]
\newtheorem{proposition}{Proposition}[section]
\newtheorem{lemma}{Lemma}[section]
\newtheorem{definition}{Definition}[section]
\newtheorem{corollary}{Corollary}[section]
\newtheorem{remark}{Remark}[section]
\newcommand{\bs}{\textcolor{black}}
\newcommand{\iii}{{\, \vert\kern-0.25ex\vert\kern-0.25ex\vert\, }}
\newcommand{\vect}[1]{\mathit{\boldsymbol{#1}}}
\newcommand{\ffi}{\varphi}
\def\Lim{\displaystyle\lim}
\def\Sup{\displaystyle\sup}
\def\Inf{\displaystyle\inf}
\def\Max{\displaystyle\max}
\def\Min{\displaystyle\min}
\def\Sum{\displaystyle\sum}
\def\Frac{\displaystyle\frac}
\def\Int{\displaystyle\int}
\def\n{|\kern -.05cm{|}\kern -.05cm{|}}
\def\Bigcap{\displaystyle\bigcap}
\def\E{{\cal E}}
\def\R{{\bf \hbox{\sc I\hskip -2pt R}}} 
\def\N{{\bf \hbox{\sc I\hskip -2pt N}}} 
\def\Z{{\bf Z}} 
\def\Q{{\bf \hbox{\sc I\hskip -7pt Q}}} 
\def\C{{\bf \hbox{\sc I\hskip -7pt C}}} 
\def\T{{\bf \hbox{\sc I\hskip -5.7pt T}}} 
{\begin{center} {\ } \vskip 1cm {\large\bf Note on Boundary Stabilization of Degenerate Schr\"{o}dinger Equations
 \\ } \

\end{center}
\begin{center}
\begin{tabular}{c}
\sc Abdelkader Benaissa and Abbes Benaissa,
\\ \\
\small  Laboratory of Analysis and Control of PDEs, \qquad \\
\small Djillali Liabes University,\qquad \\
\small P. O. Box 89, Sidi Bel Abbes 22000, ALGERIA.\qquad \\
\small E-mails: aekbenaissa@yahoo.fr\\
\small benaissa$_{-}$abbes@yahoo.com
\end{tabular}
\end{center}

\begin{abstract}
A degenerate Schr\"{o}dinger equation under fractional integral damping is considered.
Here the damping term is singular and not integrable and we consider the two cases when
damping acting on the degenerate boundary and nondegenerate boundary.
In this paper, we establish polynomial energy decay rates for
the degenerate Schr\"{o}dinger  equation by using resolvent estimates.
\end{abstract}


\section{Introduction}
The Schr\"{o}dinger-type equations are the very important models of mathematical physics, and
they have been also  intensely appeared in material science, plasma physics, quantum physics,quantum gases,
nonlinear optics, and so on.
G. Fibich {\bf\cite{fibic}} noted that the damping (absorption) term plays a significant impact in
the real models and it is more suitable to not be neglegted. This motivated several authors to deal with boundary
stabilization of Schr\"{o}dinger-type equations without degeneracies and
many decay estimates has been realized (see {\bf\cite{Ire-Trig}}).

On the contrary, when the principal part is degenerate not much is known in the documentation,
despite that many problems that are  relevant for applications are modeled by degenerate Schr\"{o}dinger
equations (see {\bf\cite{jawad}}, {\bf\cite{fekir}}, {\bf\cite{choua}}).

\noindent
In {\bf\cite{jawad}}, the authors examined the following degenerate Schr\"{o}dinger equation
\begin{equation}\label{zaz}
\left\{
\matrix{i y_t(x,t)+(x^\alpha y_x(x,t))_x=0\hfill & \hbox{ in } (0,1)\times(0,+\infty),\hfill & \cr
y(0,t)=0\hfill &  \hbox{ on } (0,+\infty),\hfill & \cr
y_t(1,t)+y_x(1,t)+y(1,t)=0\hfill & \hbox{ on } (0,+\infty),\hfill & \cr
y(x,0)=y_0(x)\hfill & \hbox{ on } (0,1).\hfill & \cr}\right.
\end{equation}
They showed that the solution decays exponentially in a suitable energy space.
Moreover, the degeneracy does not affect the decay rates of the energy.

Recently, {\bf\cite{fekir}} and {\bf\cite{choua}} explored the well-posedness and asymptotic stability of the following
degenerate Schr\"{o}dinger equation
$$
(P)\ \left\{
\matrix{y_t(x,t)-\imath(x^\alpha y_x(x,t))_x=0,\hfill & \cr
(x^\alpha y_x)(0,t)=-\imath\rho \partial^{\beta, \gamma} y(0,t),\hfill & \cr
y_x(1,t)=0,\hfill & \cr
y(x,0)=y_0(x),\hfill & \cr}\right.
(P')\ \left\{
\matrix{
y_t(x,t)-\imath(\kappa(x) y_x(x,t))_x=0, \hfill &  \hbox{ in } (0,1)\times(0,+\infty),  \hfill &\cr
\left\{\matrix{y(0, t)=0 \hfill & \hbox{ if } 0\leq m_{\kappa}< 1 \hfill & \cr
(\kappa(x) y_{x})(0,t)=0 \hfill & \hbox{ if } 1\leq m_{\kappa}< 2 \hfill & \cr}\right.\hfill  & \hbox{ in } (0, +\infty) ,\hfill&\cr
(\kappa y_x)(1,t)=\imath\rho \partial^{\beta, \gamma}y(1,t) \hfill & \hbox{ in } (0,+\infty),\hfill & \cr
y(x,0)=y_0(x) \hfill &    \hbox{ on } (0,1),\hfill & \cr}
\right.
$$
where $0<\alpha<1, \rho>0$ and $\kappa \in {\cal C}([0,1])\cap {\cal C}^1(]0,1])$ a function  satisfying the following
hypotheses:
$$
\quad \left \{
\begin{array}{ll}
\kappa(x)>0, \ \forall \ x \in ]0,1], \ \kappa(0)=0,\\
m_{\kappa}=\Sup_{x<0\leq 1} \Frac{x|\kappa'(x)|}{\kappa(x)}<2,\quad
\kappa \in {\cal C}^{[m_{\kappa}]}([0,1]).
\end{array}
\right.
$$
The term $\partial^{\beta, \gamma}$ stands for the generalized fractional integral of order $0<\beta< 1$,
which is given by
$$
\partial^{\beta, \gamma}w({t})=
\Frac{1}{\Gamma(1-\beta)}\Int_{0}^{{t}} ({t}-s)^{-\beta}e^{-{\gamma}({t}-s)}w(s)\, ds,\hfill  \hbox{ for }  \gamma\geq 0.
\label{cap}
$$
We outline the main results of {\bf\cite{fekir}} and {\bf\cite{choua}} in the following:

\noindent
The semigroup $(e^{t{\cal A}_j})_{t\geq 0},\ j=1, 2$ has the following stability properties

\noindent
(i) If $\gamma> 0$ and $\beta<(4-3\alpha)/(4-2\alpha)$, it is polynomial of order
$\frac{1}{(4-3\alpha)/(4-2\alpha)-\beta}$ (that is \\ $\|(i\lambda-{\cal A}_1)^{-1}\|_{{\cal H}}\leq c |\lambda|^{(4-3\alpha)/(4-2\alpha)-\beta}$ as $|\lambda|\rightarrow +\infty$).\\
(i) If $\gamma> 0$ and $\beta\geq(4-3\alpha)/(4-2\alpha)$, it is exponential.\\
(iii) If $\gamma> 0$, it is polynomial of order $\frac{1}{1-\beta}$ (that is $\|(i\lambda-{\cal A}_2)^{-1}\|_{{\cal H}}\leq c |\lambda|^{1-\beta}$ as $|\lambda|\rightarrow +\infty$).\\
(iv) If $\gamma=0$, it is strongly stable.

The main result of this paper then concerns the precise asymptotic behaviour of the solutions to the problems $(P)$ and $(P')$
when $\gamma=0$, with initial condition in a special subspace of $D({\cal A}_j)$.
Our technique is a special frequency and spectral analysis of the corresponding operator.
\section{Preliminaries}
Both systems $(P)$ and $(P')$ are transformed to the equivalent systems where the boundary conditions
$(x^\alpha y_x)(0,t)=-\imath\rho \partial^{\beta, 0} y(0,t)$ and
$(\kappa y_x)(1,t)=\imath\rho \partial^{\beta, 0}y(1,t)$ are replaced by the following
$$
\left\{\matrix{\psi_t(\xi, t)+\xi^{2}\psi(\xi, t) -y(0, t)\eta(\xi)=0, \hfill & \cr
(x^\alpha y_x)(0,t)=-i\zeta \Int_{-\infty}^{+\infty}\eta(\xi)\psi(\xi, t)\, d\xi\hfill & \cr
\psi(\xi, 0)=0,\hfill & \cr}\right.\ \
\left\{\matrix{\psi_t(\xi, t)+\xi^{2}\psi(\xi, t) -y(1, t)\eta(\xi)=0, \xi\in \R,\hfill & \cr
(\kappa y_x)(1,t)=i\zeta \Int_{-\infty}^{+\infty}\eta(\xi)\psi(\xi, t)\, d\xi ,\hfill & \cr
\psi(\xi, 0)=0,\hfill & \cr}\right.
$$
where $\eta(\xi)=|\xi|^{(2\beta-1)/2}$ and $\zeta=\rho (\pi)^{-1}\sin(\beta\pi)$.
The energy of solutions of $(P)$ and $(P')$ at the time $t>0$ are
defined respectively by
\begin{equation}\label{27'}
  E_j(t)=\frac{1}{2}\int_{0}^{1}|y(x,t)|^2\, dx+\frac{\zeta}{2}\Int_{-\infty}^{+\infty}|\psi(\xi, t)|^2\, d\xi\quad (j=1, 2).
\end{equation}
By Green's formula we can prove that for all $t> 0$, we have
\begin{equation}
E'_j(t)=-\zeta \Int_{-\infty}^{+\infty}\xi^{2}|\psi(\xi,t)|^2d\xi\leq 0,\quad (j=1, 2).
\label{e11}
\end{equation}
We introduce the Hilbert space
$$
{{\cal H}}=L^{2}(0,1)\times L^2(\R),
$$
with the following scalar product
$$
\langle {\cal Y},\widetilde{{\cal Y}} \rangle_{{\cal H}}=\int_{0}^{1}y(x)\overline{\widetilde{y}}(x)dx+\zeta\Int_{-\infty}^{+\infty}\psi(\xi)\overline{\widetilde{\psi}}(\xi)\, d\xi
$$
for all ${\cal Y}, \widetilde{{\cal Y}} \in {{\cal H}}$ with ${\cal Y}=(y, \psi)^{T}$ and $\widetilde{{\cal Y}}=(\widetilde{y}, \widetilde{\psi})^{T}$.
problems $(P)$ and $(P')$ can be written as
\begin{equation}\label{10}
{\cal Y}_t={\cal A}_j {\cal Y},\quad
{\cal Y}(0)={\cal Y}_0,
\end{equation}
where the operator ${{\cal A}_j}$
is defined by
$$
{\cal A}_1{\cal Y}=\pmatrix{\imath (x^\alpha y_x)_x  \cr
-\xi^2\psi+\eta(\xi)y(0) \cr}, \quad
{\cal A}_2{\cal Y}=\pmatrix{\imath (\kappa(x) y_x)_x  \cr
-\xi^2\psi+\eta(\xi)y(1) \cr}
$$
with domain
$$
\matrix{D({\cal A}_1)=\left\{\matrix{ (y, \psi)\in {\cal H}: y\in H_{\alpha}^{2}(0, 1),
 y_x(1)=0, \ (x^\alpha y_x)(0)=-\imath\zeta\Int_{-\infty}^{\infty}\eta(\xi)\psi(\xi)\, d\xi, & \cr
-\xi^2\psi+\eta(\xi)y(0)\in L^2(\R), |\xi|\psi\in L^2(\R)& \cr}\right\}, &\cr
D({\cal A}_2)=\left\{\matrix{ (y, \psi)\in {\cal H}: y\in H_{\kappa}^{2}(0, 1)\cap W_{\kappa}^1(0, 1),
(\kappa(x) y_x)(1)=\imath\zeta\Int_{-\infty}^{\infty}\eta(\xi)\psi(\xi)\, d\xi, & \cr
-\xi^2\psi+\eta(\xi)y(1)\in L^2(\R), |\xi|\psi\in L^2(\R)& \cr}\right\}, &\cr}
$$
here
$$
\matrix{W_{\kappa}^{1}(0,1)=\left\{\matrix{H^1_{0,\kappa}(0,1) \hfill & \hbox{ if }0\leq m_\kappa< 1,\hfill &\cr
H_{\kappa}^{1}(0,1)\hfill & \hbox{ if }1\leq m_\kappa< 2.\hfill &\cr}\right.
\left\{\matrix{H^1_{0,\kappa}(0,1)=\left\{y \in H^1_{\kappa}(0,1), \ y(0)=0\right\},\hfill \cr
H^2_{\kappa}(0,1)=\left\{y \in L^2(0,1), y\in H^1_{\kappa}(0,1), \ \kappa(x)y_x\in H^1(0,1)\right\}}\right.
\hfill &\cr
H^1_{\kappa}(0,1)=\left\{y \in L^2(0,1), y\hbox{ is locally absolutely continuous in }(0,1], \ \sqrt{\kappa(x)}y_x\in L^2(0,1)\right\},\hfill &\cr
H^2_{\alpha}(0,1)=\left\{y \in L^2(0,1), y\in H^1_{\alpha}(0,1), \ x^{\alpha}y_x\in H^1(0,1)\right\},\hfill &\cr
H^1_{\alpha}(0,1)=\left\{y \in L^2(0,1), y\hbox{ is \hbox{lac} in }(0,1], \ x^{\alpha/2}y_x\in L^2(0,1)\right\}.\hfill &\cr}
$$
\subsubsection*{Wellposedness and Strong Stability [{\bf\cite{fekir}} and {\bf\cite{choua}}].}
In {\bf\cite{fekir}} and {\bf\cite{choua}}, it was established that both problems $(P)$ and $(P')$ are well-posedness
and strong asymptotic stables for $\gamma=0$. The spectrum $\sigma({\cal A}_j)$ of ${\cal A}_j$
consists of isolated eigenvalues and satisfies $\sigma({\cal A}_j)\cap i\R=\{0\}$.
We apply a result due to Batty, Chill and Tomilov (\cite{batty}, Theorem 7.6) which link the decay of the associate
semigroup to the growth of $(i\lambda - {\cal A}_j)^{-1}$ near zero.
\begin{theorem}[\cite{batty}]
Let $S(t)$ be a bounded $C_0$-semigroup on a Hilbert space ${\cal X}$ with generator ${\cal A}$.
Assume that $\sigma({\cal A})\cap i\R=\{0\}$ and that there exist
$\vartheta\geq 1$ and $\upsilon > 0$ such that
$$
\|(is I-{\cal A})^{-1}\|_{{\cal L}({\cal X})}=\left\{\matrix{O(|s|^{-\vartheta}),\quad & s\rightarrow 0, \hfill & \cr
O(|s|^{\upsilon}),\quad & |s|\rightarrow \infty. \hfill & \cr}\right.
$$
Then there exist constants $C, t_0> 0$ such that for all $t\geq t_0$ and $U_0\in D({\cal A})\cap R({\cal A})$
we have
$$
\|e^{{\cal A}t}U_0\|^{2}\leq C\Frac{1}{t^{\frac{2}{\varsigma}}}
\|U_0\|_{D({\cal A)}\cap R({\cal A})}^{2}\quad (\hbox{ with }\varsigma=\max\{\vartheta, \upsilon\}).
$$
\label{thm2kkk}
\end{theorem}
\section{Main result}

Our main results are the following.
\begin{theorem}
The semigroup ${S_{{\cal A}_1}(t)}_{t\geq 0}$ is polynomially stable and
$$
E_1(t)=\|S_{{\cal A}_1}(t){\cal Y}_0\|_{\cal H}^2\leq \Frac{1}{t^{2}}  \|{\cal Y}_0\|_{D({\cal A}_1)\cap R({\cal A}_1)}^2.
$$
\label{thm225kkk}
\end{theorem}
\begin{theorem}
The semigroup ${S_{{\cal A}_2}(t)}_{t\geq 0}$ is polynomially stable and
$$
E_2(t)=\|S_{{\cal A}_2}(t){\cal Y}_0\|_{\cal H}^2\leq \Frac{1}{t^{\frac{2}{2-\beta}}}  \|{\cal Y}_0\|_{D({\cal A}_2)\cap R({\cal A}_2)}^2.
$$
Moreover in the particular case $\kappa(x)=x^{\alpha}, 0<\alpha <2$, we obtain better decay estimate that is
$$
E_2(t)=\|S_{{\cal A}_2}(t){\cal Y}_0\|_{\cal H}^2\leq \Frac{1}{t^{2}}  \|{\cal Y}_0\|_{D({\cal A}_2)\cap R({\cal A}_2)}^2.
$$
\label{thmbb}
\end{theorem}
{\bf Proof of Theorem \ref{thm225kkk}}
We will study the behaviour of the resolvent operator $(i\lambda-{{\cal A}_1})^{-1}$ as $|\lambda|\rightarrow 0$.
For any $F=(f_1, f_2)^T\in {{\cal H}}$, we examine the equation
\begin{equation}\label{56}
(\imath\lambda-{{\cal A}_1}){\cal Y}=F.
\end{equation}
Clearly, ${\cal Y}\in D({{\cal A}_1})$ is a solution of the equation if and only if
\begin{equation}\label{58}
\left\{\matrix{-\lambda y+(x^\alpha y_x)_x=\imath f_1,  \hfill & \cr
i\lambda \psi+\xi^{2}\psi-y(0)\eta(\xi)=f_2,\hfill & \cr
(x^\alpha y_x)(0)=-i\rho (i\lambda)^{\beta-1} y(0)
-i\zeta\Int_{-\infty}^{+\infty}\Frac{\eta(\xi)f_2(\xi)}{i\lambda+\xi^{2}}\, d\xi,\hfill & \cr
y_x(1)=0.\hfill & \cr
}\right.
\end{equation}
The general solution of the differential equation $(\ref{58})_1$ can be written as
\begin{equation}
y(x)= A\theta_+(x)+B\theta_-(x)-\frac{\pi}{2\sin \nu_\alpha \pi}
\left(\frac{2}{2-\alpha}\right)\int_{0}^{x}\imath f_1(X)(\theta_+(X)
\theta_-(x)-\theta_+(x)\theta_-(X))dX,
\label{eqtt}
\end{equation}
where $\theta_+$ and $\theta_-$ are defined by
\begin{equation}\label{theta}
\theta_+(x)=x^{\frac{1-\alpha}{2}}J_{\nu_{\alpha}}\left(\frac{2}{2-\alpha}\mu x^{\frac{2-\alpha}{2}}\right)
\quad \hbox{ and  } \
\theta_-(x)=x^{\frac{1-\alpha}{2}}J_{-\nu_{\alpha}}\left(\frac{2}{2-\alpha}\mu x^{\frac{2-\alpha}{2}}\right)
\end{equation}
with $\mu=i\sqrt{\lambda}$, $\nu_\alpha=\frac{1-\alpha}{2-\alpha}$. $J_{\nu_{\alpha}}$ and $J_{-\nu_{\alpha}}$ are Bessel functions of the first kind of order $\nu_{\alpha}$ and $-\nu_{\alpha}$.
Using $(\ref{58})_3$ and $(\ref{58})_4$, the constants $A$ and $B$ satisfy
\begin{equation}\label{77}
\left( \begin{array}{c}
A \\\\ B
\end{array} \right)=\frac{1}{D}
\left( \begin{array}{cc}
\theta'_-(1) &
-\imath\rho(i\lambda)^{\beta-1} d^- \\\\
-\theta'_+(1) &(1-\alpha)d^+
\end{array}\right)
\left( \begin{array}{c}
C \\\\ \tilde{C}
\end{array} \right),
\end{equation}
where
$$
\theta_+'(1)=(1-\alpha)J_{\nu_{\alpha}}\left(\frac{2\mu}{2-\alpha}\right)-\mu J_{\nu_{\alpha}+1}\left(\frac{2\mu}{2-\alpha}\right),\quad
\theta_-'(1)=-\mu J_{-\nu_{\alpha}+1}\left(\frac{2\mu}{2-\alpha}\right),
$$
\begin{equation}\label{dd}
d^+=c^+_{\nu_{\alpha,0}}\left(\frac{2}{2-\alpha}\mu\right)^{\nu_{\alpha}} \ \hbox{ and \ \ }
d^-=c^-_{\nu_{\alpha,0}}\left(\frac{2}{2-\alpha}\mu\right)^{{-\nu_{\alpha}}}
\end{equation}
and
$$
C=-i\zeta\Int_{-\infty}^{+\infty}\Frac{\eta(\xi)f_2(\xi)}{i\lambda+\xi^{2}}\, d\xi,\quad
\tilde{C}=\frac{\pi}{2\sin \nu_\alpha \pi}
\left(\frac{2}{2-\alpha}\right)\int_{0}^{1}\imath f_1(X)(\theta_+(X)
\theta'_-(1)-\theta'_+(1)\theta_-(X))dX.
$$
Solving $(\ref{77})$ finally gives
$$
|A|=\left| \frac{\theta_-'(1)C-\imath\rho(i\lambda)^{\beta-1} \tilde{C}d^-}{D}\right|
\hbox{ and }
|B|=\left| \frac{-\theta_+'(1)C+(1-\alpha)\tilde{C}d^+}{D}\right|,
$$
where
$$
D=(1-\alpha)d^+\theta'_-(1)-\imath\rho\theta'_+(1)(i\lambda)^{\beta-1} d^-.
$$
Then, using Cauchy-Schwarz inequality, the expressions of $\theta_+'$ and $\theta_-'$, we get
\begin{equation}\label{qqa}
\matrix{|C| &\leq & \zeta \left(\Int_{-\infty}^{\infty}\Frac{\eta(\xi)^2}{|i\lambda+\xi^2|^2}\, d\xi\right)^{1/2}\|f_2\|_{L^2(-\infty, +\infty)}\hfill \cr
&\leq & \sqrt{2}\zeta \left(\Int_{-\infty}^{\infty}\Frac{\eta(\xi)^2}{(|\lambda|+\xi^2)^2}\, d\xi\right)^{1/2}\|f_2\|_{L^2(-\infty, +\infty)}
\leq  c|\mu|^{\beta-2}\|f_2\|_{L^2(-\infty, +\infty)},\hfill \cr}
\end{equation}
\begin{equation}\label{C}
|\tilde{C}|\leq c\|f_1\|_{L^2(0,1)}.
\end{equation}
Indeed, in the neighborhood of $\mu=0$ we have
$$
\theta_+'(1)\sim (1-\alpha)c^+_{\nu_{\alpha,0}}\left(\frac{2}{2-\alpha}\mu\right)^{\nu_{\alpha}}, \quad
\theta_-'(1)\sim -\mu c^-_{\nu_{\alpha,0}}\left(\frac{2}{2-\alpha}\mu\right)^{{-\nu_{\alpha}}},\ \ \ \
|D|\geq C|\mu|^{2\beta-2}.
$$
and
$$
\|\theta_{+}\|_{L^{2}(0, 1)}^2=\Frac{1}{2-\alpha}\Frac{1}{r^2}\left[\left(r J_{\nu_{\alpha}}(r)\right)^2
+\left(r J_{\nu_{\alpha}+1}(r)\right)^2
-2\nu_{\alpha}r J_{\nu_{\alpha}}\left(r\right)J_{\nu_{\alpha}+1}\left(r\right)\right]
\sim \Frac{1}{2-\alpha}(c_{\nu_{\alpha}, 0}^+)^2 r^{2\nu_{\alpha}},
$$
where $r=\frac{2\mu}{2-\alpha}$. Using $(\ref{qqa})$, $(\ref{C})$ and $(\ref{dd})$, we deduce that
$$
|A|\leq c'|\mu|^{-\nu_{\alpha}}(\|f_1\|_{L^2(0,1)}+\|f_2\|_{L^2(-\infty, +\infty)})
$$
and
\begin{equation}\label{qe8}
|B|\leq  c |\mu|^{\nu_{\alpha}-\beta}(\|f_1\|_{L^2(0,1)}+\|f_2\|_{L^2(-\infty, +\infty)}).
\end{equation}
Then
\begin{equation}\label{qe12}
\|y\|_{L^2(0,1)}\leq c |\mu|^{-\beta}(\|f_1\|_{L^2(0,1)}+\|f_2\|_{L^2(-\infty, +\infty)}).
\end{equation}
Moreover from $(\ref{58})_2$, we have $\psi=\Frac{y(0)\eta(\xi)+f_2(\xi)}{i\lambda+\xi^2}$.\\
Then
$$
\matrix{\|\psi\|_{L^2(-\infty, +\infty)}^2&\leq & 2|y(0)|^2\Int_{-\infty}^{+\infty}\Frac{\eta(\xi)^2}{|i\lambda+\xi^2|^2}\, d\xi+2\Int_{-\infty}^{+\infty}\Frac{|f_2(\xi)|^2}{|i\lambda+\xi^2|^2}\, d\xi\hfill \cr
&\leq & 2 \Frac{\pi}{\sin\frac{\beta}{2}\pi} |\lambda|^{\beta-2}|y(0)|^2+\Frac{1}{|\lambda|^2}\|f_2\|_{L^2(-\infty, +\infty)}^{2}.\hfill \cr}
$$
From (\ref{eqtt}) and (\ref{qe8}), it is clear that
$|y(0)|^2= |B|^2 (d^{-})^2\leq c|\mu|^{-2\beta}(\|f_1\|_{L^2(0,1)}^2+\|f_2\|_{L^2(-\infty, +\infty)}^2)$.\\
Hence
\begin{equation}\label{qe13}
\|\psi\|_{L^2(-\infty, +\infty)}^2\leq \Frac{c}{|\lambda|^2} (\|f_1\|_{L^2(0,1)}^2+\|f_2\|_{L^2(-\infty, +\infty)}^2).
\end{equation}
We deduce by (\ref{qe12}) and (\ref{qe13}) that
$$
\|{\cal Y}\|_{{\cal H}}^2=\|y\|_{L^2(0,1)}^2+\|\psi\|_{L^{2}(-\infty, \infty)}^2\leq  c'\left(\Frac{1}{|\lambda|^{\beta}}+\Frac{1}{|\lambda|^2}\right) \|F\|_{\cal H}^2
\leq  \Frac{c}{|\lambda|^2} \|F\|_{\cal H}^2.
$$
Finally, we obtain $\|(i\lambda -{\cal A}_1)^{-1}\|_{{\cal H}}\leq \frac{c}{|\lambda|}$ as $\lambda\rightarrow 0$.

\noindent
{\bf Proof of Theorem \ref{thmbb}} We will need to study the
resolvent equation $(i\lambda-{\cal A}_2){\cal Y}=F$, for $\lambda\in \R$,
namely
\begin{equation}
\left\{\matrix{i\lambda y-\imath (\kappa(x) y_x)_x=f_1,  \hfill & \cr
i\lambda \psi+\xi^{2}\psi-y(1)\eta(\xi)=f_2.\hfill & \cr
\matrix{
\left\{\matrix{y(0, t)=0 \hfill & \hbox{ if } 0\leq m_{\kappa}< 1 \hfill & \cr
(\kappa(x)y_{x})(0,t)=0 \hfill & \hbox{ if } 1\leq m_{\kappa}< 2 \hfill & \cr}\right.
\hfill  & \hbox{ in } (0, +\infty) ,\hfill&\cr
(\kappa y_x)(1,t)-\imath\zeta \Int_{-\infty}^{+\infty}\eta(\xi)\psi(\xi, t)\, d\xi=0 \hfill & \hbox{ in } (0,+\infty),\hfill & \cr}
}\right.
\label{e188}
\end{equation}
Using $(\ref{e188})_2$ and repeating calculations in {\bf\cite{choua}} near $\lambda=0$, we get
\begin{equation}
|y(1)|^2\leq c|\lambda|^{1-\beta}\|{\cal Y}\|_{\cal H}\|F\|_{\cal H}+\Frac{c}{|\lambda|^{\beta}}\|F\|_{\cal H}^{2}.
\label{e41nkk}
\end{equation}
\begin{equation}
\left|\Int_{-\infty}^{+\infty}\eta(\xi)\psi(\xi)\, d\xi\right|^2\leq c|\lambda|^{\beta-1}\|{\cal Y}\|_{\cal H}\|F\|_{\cal H}+c|\lambda|^{\beta-2}\|F\|_{\cal H}^{2}.
\label{qe9}
\end{equation}
Let us multiply the equation $(\ref{e188})_1$ by $\overline{y}$ integrating over $(0,1)$ and integrating par parts we get
\begin{equation}\label{ma}
\lambda\|y\|^2_{L^2(0,1)}-\left[\kappa(x)y_x\overline{y}\right]_0^1+\int_{0}^{1}\kappa(x)|y_x|^2dx
=-\imath\int_{0}^{1}f_1\overline{y}dx.
\end{equation}
If $\lambda> 0$ near $0$, we have from (\ref{ma}), (\ref{e41nkk}) and (\ref{qe9}) that
$$
\matrix{\lambda\|y\|^2_{L^2(0,1)}\leq  \|f_1\|_{L^2(0, 1)}\|y\|_{L^2(0, 1)}+|(\kappa(x)y_x)(1)||\overline{y}(1)|\hfill \cr
\leq  \|f_1\|_{L^2(0, 1)}\|y\|_{L^2(0, 1)}+\zeta|\Int_{\R}\eta(\xi)\psi(\xi)\, d\xi||\overline{y}(1)|\hfill \cr
\leq  \|F\|_{\cal H}\|{\cal Y}\|_{\cal H}+(c|\lambda|^{\frac{\beta-1}{2}}\|F\|_{\cal H}\|{\cal Y}\|_{\cal H})^{1/2}+c|\lambda|^{\frac{\beta-2}{2}}\|F\|_{\cal H})(|\lambda|^{\frac{1-\beta}{2}}(\|{\cal Y}\|_{\cal H}\|F\|_{\cal H})^{1/2}+c|\lambda|^{-\frac{\beta}{2}}\|F\|_{\cal H}).\hfill \cr
\leq  c\|F\|_{\cal H}\|{\cal Y}\|_{\cal H}+c|\lambda|^{-\frac{1}{2}}(\|F\|_{\cal H}\|{\cal Y}\|_{\cal H})^{1/2}|\|F\|_{\cal H}+\Frac{c}{|\lambda|}\|F\|_{\cal H}^2
\leq  c\|F\|_{\cal H}\|{\cal Y}\|_{\cal H}+\Frac{c}{|\lambda|}\|F\|_{\cal H}^2.\hfill \cr
}
$$
Moreover
\begin{equation}
\matrix{\|\psi\|_{L^2(\R)}^2&\leq & 2|y(1)|^2\Int_{\R}\Frac{\eta(\xi)^2}{|i\lambda+\xi^2|^2}\, d\xi+2\Int_{\R}\Frac{|f_2(\xi)|^2}{|i\lambda+\xi^2|^2}\, d\xi\hfill \cr
&\leq & 2 \Frac{\pi}{\sin\frac{\beta}{2}\pi} |\lambda|^{\beta-2}|y(1)|^2+\Frac{1}{|\lambda|^2}\|f_2\|_{L^2(\R)}^{2}
\leq   c |\lambda|^{-1}\|F\|_{\cal H}\|{\cal Y}\|_{\cal H}+ \Frac{c}{|\lambda|^2}\|F\|_{\cal H}^2\hfill \cr}
\label{qe1}
\end{equation}
Then, we deduce that $\|{\cal Y}\|_{\cal H}^2\leq \frac{C}{\lambda^2} \|F\|_{\cal H}^2$.
Now, if $\lambda< 0$ near $0$,
we multiply (\ref{e188}) by $-2x\overline{y}_x$, integrating over $(0,1)$, and integrating par parts we obtain
$$
\matrix{-\lambda\left[x|y|^2\right ]_0^1+\lambda\|y\|^2_{L^2(0,1)}+2\left[\kappa(x)|y_x|^2x\right]_0^1-2\int_{0}^{1}\kappa(x)|y_x|^2dx-\left[x\kappa(x)|y_x|^2\right]_0^1\hfill &\cr
+\int_{0}^{1}\kappa(x)|y_x|^2dx+\int_{0}^{1}x\kappa'(x)|y_x|^2dx=2\Re\ i\int_{0}^{1}f_1 x\overline{y}_xdx.\hfill &\cr}
$$
Multiplying $(\ref{ma})$ by $-\frac{m_{\kappa}}{2}$, summing with the last equation and taking the real part, we get
$$
\matrix{
\lambda\left(1-\frac{m_{\kappa}}{2}\right)\|y\|^2_{L^2(0,1)}-\int_{0}^{1}\left(\kappa(x)-x\kappa'(x)
+\frac{m_{\kappa}}{2}\kappa(x)\right)|y_x|^2dx-\lambda\left[x|y|^2\right]_0^1+\left[x\kappa(x)|y_x|^2\right]_0^1\hfill &\cr
+\frac{m_{\kappa}}{2}\left[\kappa(x)y_x\overline{v}\right]_0^1=2\imath\Re\int_0^1f_1 x\overline{y}_x\, dx
+\frac{m_{\kappa}}{2}\imath\int_0^1f_1\overline{y}dx.\hfill &\cr}
$$
Then, there exist $c_1$, $c_2$ and $c_3$ such that
$$
\matrix{|\lambda|\left(1-\frac{m_{\kappa}}{2}\right)\|y\|^2_{L^2(0,1)}+\int_{0}^{1}\left(\kappa(x)-x\kappa'(x)+\frac{m_{\kappa}}{2}\kappa(x)\right)|y_x|^2dx \leq \hfill &\cr
c|\lambda|(|\lambda|^{1-\beta}\|{\cal Y}\|_{\cal H}\|F\|_{\cal H}+\Frac{c}{|\lambda|^{\beta}}\|F\|_{\cal H}^{2})
+c|\lambda|^{\beta-1}\|{\cal Y}\|_{\cal H}\|F\|_{\cal H}+c|\lambda|^{\beta-2}\|F\|_{\cal H}^{2}\hfill &\cr
+(c|\lambda|^{\frac{\beta-1}{2}}\|F\|_{\cal H}\|{\cal Y}\|_{\cal H})^{1/2}+c|\lambda|^{\frac{\beta-2}{2}}\|F\|_{\cal H})(|\lambda|^{\frac{1-\beta}{2}}(\|{\cal Y}\|_{\cal H}\|F\|_{\cal H})^{1/2}+c|\lambda|^{-\frac{\beta}{2}}\|F\|_{\cal H})\hfill &\cr
+c\|{\cal Y}\|_{\cal H}\|F\|_{\cal H}+c\|F\|_{\cal H}\|\sqrt{\kappa(x)} y_x\|_{L^2(0, 1)}\hfill &\cr
\leq \Frac{c_1}{|\lambda|^{1-\beta}}\|{\cal Y}\|_{\cal H}\|F\|_{\cal H}+\Frac{c_2}{|\lambda|^{2-\beta}} \|F\|_{\cal H}^2+c_3\|F\|_{\cal H}\|\sqrt{\kappa(x)} y_x\|_{L^2(0, 1)}.\hfill &\cr
}
$$
From this, we obtain
\begin{equation}
\|y\|^2_{L^2(0,1)}\leq
\Frac{C}{|\lambda|^{2-\beta}}\|{\cal Y}\|_{\cal H}\|F\|_{\cal H}+\Frac{C}{|\lambda|^{3-\beta}} \|F\|_{\cal H}^2.
\label{qe2}
\end{equation}
Then by (\ref{qe1}) and (\ref{qe2}) we deduce that $\|{\cal Y}\|_{\cal H}\leq \frac{C}{|\lambda|^{2-\beta}}\|F\|_{\cal H}$.

\noindent
{\bf Case $\kappa(x)=x^{\alpha},\ 0<\alpha< 2$}. We treat only the case $0< \alpha <1$. The case $1\leq \alpha< 2$
is similar. Using (\ref{eqtt}) and the boundary conditions, we have $y(0)=0\Rightarrow B d^{-}=0\Rightarrow B=0$ and
$$
\matrix{A[\theta'_+(1)-i\rho (i\lambda)^{\beta-1}\theta_+(1)]=
\frac{\pi}{2\sin \nu_\alpha \pi}
\left(\frac{2}{2-\alpha}\right)\int_{0}^{1}\imath f_1(X)(\theta_+(X)
\theta'_-(1)-\theta'_+(1)\theta_-(X))dX&\hfill \cr
-i\rho (i\lambda)^{\beta-1}\frac{\pi}{2\sin \nu_\alpha \pi}
\left(\frac{2}{2-\alpha}\right)\int_{0}^{1}\imath f_1(X)(\theta_+(X)
\theta_-(1)-\theta_+(1)\theta_-(X))dX+
i\zeta\Int_{-\infty}^{+\infty}\Frac{\eta(\xi)f_2(\xi)}{i\lambda+\xi^{2}}\, d\xi&\hfill \cr}
$$
We have easily
$|D|\simeq c|\mu|^{2\beta+\nu_{\alpha}-2}$. Then
$|A|\leq c |\mu|^{-\beta-\nu_{\alpha}}(\|f_1\|_{L^2(0,1)}+\|f_2\|_{L^2(\R)})$.\\
Then
$$
\|y\|_{L^2(0,1)}\leq c |\mu|^{-\beta}(\|f_1\|_{L^2(0,1)}+\|f_2\|_{L^2(\R)}).
$$
From (\ref{eqtt}), we have
$$
|y(1)|^2= \leq c|\mu|^{-2\beta}(\|f_1\|_{L^2(0,1)}^2+\|f_2\|_{L^2(\R)}^2).
$$
Hence
$$
\|\psi\|_{L^2(\R)}^2\leq \Frac{c}{|\lambda|^2} (\|f_1\|_{L^2(0,1)}^2+\|f_2\|_{L^2(\R)}^2).
$$
Thus, we conclude that
$$
\|{\cal Y}\|_{{\cal H}}^2=\|y\|_{L^2(0,1)}^2+\|\psi\|_{L^{2}(-\infty,
\infty)}^2 \leq
c'(\Frac{1}{|\lambda|^{\beta}}+\Frac{1}{|\lambda|^2})
\|F\|_{\cal H}^2\leq  \Frac{c}{|\lambda|^2} \|F\|_{\cal H}^2.
$$
Lastly,, we obtain $\|(i\lambda -{\cal A}_2)^{-1}\|_{{\cal H}}\leq \frac{c}{|\lambda|}$ as $\lambda\rightarrow 0$.
The conclusion follows by applying Theorem $\ref{thm2kkk}$ with $\varsigma=\max\{1, 1-\beta\}=1$.

\end{document}